\documentclass{article}
\usepackage{amsmath, amssymb, amsthm}

\theoremstyle{definition}

\title{Triangular hyperbolic buildings.}

\author{Riikka Kangaslampi, Alina Vdovina}

\begin{document}

\maketitle

\abstract{ We construct triangular hyperbolic
polyhedra whose links are generalized $4$-gons.
The universal cover of those polyhedra are hyperbolic
buildings, which appartments are hyperbolic planes tesselated
by regular triangles with angles $\pi/4$. Moreover, the
fundamental groups of the polyhedra acts simply transitively
on vertices of the buildings.}

\bigskip

\section*{Introduction}

A {\em polyhedron} is a
two-dimensional complex which is obtained from several oriented
$p$-gons by identification of corresponding sides.
We will consider euclidean and hyperbolic polyhedra. Let's take a
point of the polyhedron and take a sphere (euclidean or hyperbolic)
 of
a small radius at this point. The intersection of the sphere with
the polyhedron is a graph, which is called the {\em link} at this
point.

We consider  polyhedra such that all links of all
vertices are generalized 4-gons.

We  say that a polyhedron $P$ is an
$(m,n)$-polyhedron if the girth of any link of $P$ is at
least $m$ and each face of $P$ is a polygon with at least $n$
edges.

Let $P$ be a {$(m,n)$-polyhedron such
that $m$ and $n$} satisfy the inequality $mn \geq 2(m+n)$.
 This inequality appears, for example,
 in small cancellation theory
\cite{LS}.
The universal covering of a $(m,n)$-polyhedron with the metric
introduced in \cite[p.~165]{BBr} is a complete metric space
of non-positive curvature in the sense of Alexandrov and
Busemann \cite{GH}.

In this note we construct examples of
$(4,3)$-polyhedra. It follows from \cite{BS}, that
the fundamental groups of our polyhedra satisfy
the property (T) of Kazhdan. (Another relevant reference
is \cite{Z})

\bigbreak

\noindent{\bf Definition.}
A {\em generalized $m$-gon} is a graph which is a union of subgraphs,
called apartments, such that:
\begin{itemize}
\item[1.] Every apartment is a cycle composed of $2m$ edges
for some fixed $m$.
\item[2.] For any two edges there is an apartment containing both of them.
\item[3.] If two apartments have a common edge, then there is
an isomorphism between them fixing their intersection pointwise.
\end{itemize}

\bigbreak

\noindent
{\bf Definition.} Let $\mathcal{P}(p,m)$ be a tessellation of the
hyperbolic plane by  regular polygons with $p$ sides,
with angles $\pi/m$ in each vertex where $m$ is an integer.
A {\em hyperbolic building} is a polygonal complex $X$,
which can be expressed as the union of subcomplexes called apartments
such that:

\begin{itemize}
\item[1.] Every apartment is isomorphic to $\mathcal{P}(p,m)$.
\item[2.] For any two polygons of $X$, there is an apartment
containing both of them.
\item[3.] For any two apartments $A_1, A_2 \in X$ containing
the same polygon, there exists an isomorphism $ A_1 \to A_2$
fixing $A_1 \cap A_2$.
\end{itemize}

\bigbreak

Let $C_p$ be a polyhedron whose faces are $p$-gons
and whose links are generalized $m$-gons with $mp>2m+p$.
We equip every face of $C_p$ with the hyperbolic
metric such that all sides of the polygons are geodesics
and all angles are $\pi/m$. Then the universal covering of such
polyhedron is a hyperbolic building, see \cite{Paulin}.
The polyhedra we construct in this note are from this class,
since $p=3$, $m=4$.
Our construction gives examples of hyperbolic triangular buildings
with regular triangles as chambers, which were known before,
on our knowledge. Examples of hyperbolic buildings with right-angled
triangles were constructed by M.Bourdon \cite{Bourdon}.

In the case $p=3$, $m=3$, i.e. $C_p$ is a simplex, we can give
a euclidean metric to every face. In this metric
all sides of the triangles
are geodesics. The universal coverings of these polyhedra
with the euclidean metric are euclidean buildings, see \cite{BB}, \cite{Ba}.

Recall (cf. \cite{Ronan}) that
a {\em generalized m-gon} is a connected, bipartite graph of
diameter $m$ and girth $2m$, in which each vertex lies on at least
two edges. A graph is {\em bipartite} if its set of vertices can
be partitioned into two disjoint subsets such that no two vertices
in the same subset lie on a common edge. The vertices of
 one subset we will call black vertices and the
vertices of the other subset the white ones. The {\em diameter} is
the maximum distance between two vertices and the {\em girth}  is
the length of a shortest circuit.

Let $G$ be a connected bipartite graph on $q+r$ vertices,
$q$ black vertices and $r$ white ones. Let $\mathcal{A}$ and
$\mathcal{B}$ be two alphabets on $q$ and $r$ letters respectively,
$\mathcal{A}=\{x_1, x_2, \ldots, x_q\}$ and $\mathcal{B}=\{y_1, y_2, \ldots,
y_r\}$.
We mark every black vertex with an element from
$\mathcal{A}$ and every white vertex with an element from $\mathcal{B}$.


\bigbreak

\section{Polygonal presentation.}

We recall a definition of polygonal presentation
introduced in \cite{V}.

\medbreak

\noindent {\bf Definition.} Suppose we have $n$ disjoint connected
bipartite graphs\linebreak
 $G_1, G_2, \ldots G_n$.
Let $P_i$ and $Q_i$ be the sets of black and white vertices
respectively in $G_i$, $i=1,\dots,n$; let $P=\bigcup P_i$,
$Q=\bigcup Q_i$, $P_i \cap P_j = \emptyset$,
 $Q_i \cap Q_j = \emptyset$
for $i \neq j$ and
let $\lambda$ be a  bijection $\lambda: P\to Q$.

A set $\mathcal{K}$ of $k$-tuples $(x_1,x_2, \ldots, x_k)$, $x_i \in P$,
will be called a {\em polygonal presentation} over $P$ compatible
with $\lambda$ if

\begin{itemize}

\item[(1)] $(x_1,x_2,x_3, \ldots ,x_k) \in \mathcal{K}$ implies that
   $(x_2,x_3,\ldots,x_k,x_1) \in \mathcal{K}$;

\item[(2)] given $x_1,x_2 \in P$, then $(x_1,x_2,x_3, \ldots,x_k) \in \mathcal{K}$
for some $x_3,\ldots,x_k$ if and only if $x_2$ and $\lambda(x_1)$
are incident in some $G_i$;

\item[(3)] given $x_1,x_2 \in P$, then  $(x_1,x_2,x_3, \ldots ,x_k) \in \mathcal{K}$
    for at most one $x_3 \in P$.

\end{itemize}

If there exists such $\mathcal{K}$, we will call $\lambda$ a {\em basic bijection}.

Polygonal presentations for $n=1$, $k=3$ were listed in $\cite{Cart}$
with the incidence graph of the finite projective plane of 
order two or three as the graph $G_1$.
 
Some polygonal presentations for $k > 3$ were constructed in \cite{V}.

\medskip

Now we construct  polygonal
presentations with $k=3$,$n=1$,
but the graph $G_1$ is a generalized 4-gon.

$T_1$:

$(x_1,x_2,x_7)$
  
$(x_1,x_8, x_{11})$

$(x_1,x_{14},x_5)$

$(x_2,x_4,x_{13})$

$(x_{12},x_4,x_2)$

$(x_4,x_9,x_3)$

$(x_6,x_8,x_3)$

$(x_{14},x_6,x_3)$

$(x_{12},x_{10},x_5)$

$(x_{13},x_{15},x_5)$

$(x_{12},x_9,x_6)$

$(x_{11},x_{10},x_7)$

$(x_{14},x_{13},x_7)$

$(x_9,x_{15},x_8)$

$(x_{11},x_{15},x_{10})$

$T_2$:

$(x_1,x_{10},x_1)$ 

$(x_1,x_{15},x_2)$

$(x_2,x_{11},x_9)$

$(x_2,x_{14},x_3)$

$(x_3,x_7,x_4)$

$(x_3,x_{15},x_{13})$

$(x_4,x_8,x_6)$

$(x_{12},x_{11},x_4)$

$(x_5,x_8,x_5)$

$(x_5,x_{10},x_{12})$

$(x_6,x_{14},x_6)$

$(x_7,x_{12},x_7)$

$(x_{13},x_9,x_8)$

$(x_{14},x_{15},x_9)$

$(x_{13},x_{11},x_{10})$

\medskip

Let's show, that these sets are desired polygonal 
presentations. Remark, that the smallest generalized
4-gon can be presented in the following way:
its ``points'' are pairs $(ij)$, where $i,j=1,...,6$, $i \neq j$
and ``lines'' are triples of those pairs where all $i,j$ are
different. We mark pairs $(ij)$, where $i,j=1,...,6$, $i \neq j$
by numbers from 1 to 15 in natural order. Now one can
check by direct examination, that the graph $G_1$
is really the smallest generalized 4-gon. 
The polygonal presentation $T_1$ and $T_2$ are not equivalent
since there is no automorphism of the generalized
4-gon, which transforms one to another.

\section{Construction of polyhedra.}

We can associate  a polyhedron $K$ on $n$ vertices with
each polygonal presentation $\mathcal{K}$ as follows:
for every cyclic $k$-tuple $(x_1,x_2,x_3,\ldots,x_k)$ from
the definition
we take an oriented $k$-gon on the boundary of which
the word $x_1 x_2 x_3\ldots x_k$ is written. To obtain
the polyhedron we identify the corresponding sides of our
polygons, respecting orientation.
We will say that the
polyhedron $K$ {\em corresponds} to the polygonal
presentation $\mathcal{K}$.
\medskip

\noindent {\bf Lemma \cite{V}} A polyhedron $K$ which
corresponds to a polygonal presentation $\mathcal{K}$ has
  graphs $G_1, G_2, \ldots, G_n$ as the links.

\medskip

\noindent
{\bf Remark.} Consider a  polygonal
presentation $\mathcal{K}$. Let $s_i$ be the number of vertices
of the graph $G_i$ and $t_i$ be the number of edges of $G_i$,
$i=1,\dots,n$.
If the polyhedron $K$  corresponds to the polygonal
presentation $\mathcal{K}$, then $K$ has $n$ vertices
(the number of vertices of $K$ is equal to the number of graphs),
$k \sum_{i=1}^n s_i$ edges and $\sum_{i=1}^n t_i$
faces, all faces are polygons with $k$ sides.

\medskip

For polygonal presentation $T_i,i=1,2$
 take 15 oriented regular hyperbolic triangles
with angles $\pi/4$, whrite words from the presentation
on their boundaries and glue together respecting
orientation sides with the same letters. 
The result is a hyperbolic polyhedron with one vertex and
15 faces and the universal covering is a triangular
hyperbolic building. The fundamental group of the polyhedron
acts simply transitively on vertices of the building.
The group has 15 generators and 15 relations, which
come naturally from the polygonal presentation.

\end{document}